\magnification=\magstep1
\tolerance=2000

\baselineskip=17pt 

\font\tenopen=msbm10
\font\sevenopen=msbm7
\font\fiveopen=msbm5
\newfam\openfam
\def\openo{\fam\openfam\tenopen}
\textfont\openfam=\tenopen
\scriptfont\openfam=\sevenopen
\scriptscriptfont\openfam=\fiveopen

\font\tenujsym=msam10
\font\sevenujsym=msam7
\font\fiveujsym=msam5
\newfam\ujsymfam

\textfont\ujsymfam=\tenujsym
\scriptfont\ujsymfam=\sevenujsym
\scriptscriptfont\ujsymfam=\fiveujsym

\def\R{{\openo R}}
\def\N{{\openo N}}
\def\C{{\openo C}}

\def\loc{{\rm loc}} 

\def\e{\varepsilon} 

\def\kocka{\kern+0.9pt{\sqcup\kern-6.5pt\sqcap}\kern+0.9pt}

\centerline{\bf On the regular convergence of multiple series of numbers and} 

\centerline{\bf multiple integrals of locally integrable functions over $\overline{\R}^m_+$} 

\bigskip

\centerline{FERENC M\'ORICZ} 

\bigskip

University of Szeged, Bolyai Institute, Aradi v\'ertan\'uk tere 1, 

Szeged 6720, Hungary 

e-mail: moricz@math.u-szeged.hu

\vglue1cm

\noindent {\bf Abstract.} We investigate the regular convergence of the $m$-multiple series 
$$\sum^\infty_{j_1=0} \sum^\infty_{j_2=0} \ldots \sum^\infty_{j_m=0} \ c_{j_1, j_2, \ldots, j_m}\leqno(*)$$
of complex numbers, where $m\ge 2$ is a fixed integer. We prove Fubini's theorem in the discrete setting as follows. If the 
multiple series ($*$) converges regularly, then its sum in Pringsheim's sense can be 
computed by successive summation. 

We introduce and investigate the regular convergence of the $m$-multiple integral 
$$\int^\infty_0 \int^\infty_0 \ldots \int^\infty_0 f(t_1, t_2, \ldots, t_m) dt_1 dt_2 \ldots dt_m,\leqno(**)$$
where $f: \overline{\R}^m_+ \to \C$ is a locally integrable function in Lebesgue's sense over the closed positive octant 
$\overline{\R}^m_+:= [0, \infty)^m$. Our main result is a generalized version of Fubini's theorem on 
successive integration formulated in Theorem 4.1 as follows. If $f\in L^1_{\loc} (\overline{\R}^m_+)$, 
the multiple integral ($**$) converges regularly, and $m=p+q$, where $m, p\in \N_+$, then the finite limit 
$$\lim_{v_{p+1}, \ldots, v_m \to \infty} \int^{v_1}_{u_1} \int^{v_2}_{u_2} \ldots \int^{v_p}_{u_p} \int^{v_{p+1}}_0 
\ldots \int^{v_m}_0 f(t_1, t_2, \ldots, t_m) dt_1 dt_2\ldots dt_m$$
$$=:J(u_1, v_1; u_2, v_2; \ldots; u_p, v_p), \quad 0\le u_k\le v_k<\infty, \ k=1,2,\ldots, p,$$
exists uniformly in each of its variables, and the finite limit 
$$\lim_{v_1, v_2, \ldots, v_p\to \infty} J(0, v_1; 0, v_2; \ldots; 0, v_p)=I$$
also exists, where $I$ is the limit of the multiple integral ($**$) in Pringsheim's sense. 

\bigskip


\bigskip


\vglue1cm 

\noindent {\bf 1. Convergence of multiple series of numbers} 

We consider the $m$-{\it multiple series} 
$$\sum^\infty_{j_1=0} \ \sum^\infty_{j_2=0} \ldots \sum^\infty_{j_m=0} c_{j_1, j_2, \ldots, j_m} \leqno(1.1)$$
of complex numbers, where $m\ge 2$ is a fixed integer. The $m$-fold {\it rectangular partial sums} 
of (1.1) are defined by 
$$s(l_1, l_2, \ldots, l_m):= \sum^{l_1}_{j_1=0} \sum^{l_2}_{j_2=0} \ldots \sum^{l_m}_{j_m=0} c_{j_1, j_2, \ldots, j_m}, 
\quad (l_1, l_2, \ldots, l_m)\in \N^m.$$
We recall that the multiple series (1.1) is said to {\it converge in Pringsheim's sense} to the sum 
$s\in \C$, in symbols: 
$$\lim_{l_1, l_2, \ldots, l_m\to \infty} s(l_1, l_2, \ldots, l_m) = s,\leqno(1.2)$$
if for every $\e>0$ there exists $\lambda_1=\lambda_1(\e)\in \N$ such that 
$$|s(l_1, l_2, \ldots, l_m) - s| < \e \quad {\rm if} \quad 
\min\{l_1, l_2, \ldots, l_m\}>\lambda_1.\leqno(1.3)$$
See [9,10] by Pringsheim, and see also [13 on p. 303, just after formula (1.18)] by Zygmund, without indication of the term 
``in Pringsheim's sense". 

\bigskip

\noindent {\bf Remark 1.1.} In contrast to convergent single series, the convergence of multiple series in 
Pringsheim's sense implies neither the boundedness of its terms $c_{j_1, j_2, \ldots, j_m}$, nor the convergence of its subseries. 
We recall that for any choice $1\le p_1 < p_2 < \ldots < p_e \le m$ of integers, where 
$1\le e<m$, while denoting by $1\le p_{e+1} < \ldots < p_m \le m$ the remaining integers between 1 
and $m$, and for any fixed values of $(j_{p_{e+1}}, \ldots, j_{p_m}) \in \N^{m-e}$, the $e$-multiple series 
$$\sum^\infty_{j_{p_1}=0} \sum^\infty_{j_{p_2}=0} \ldots \sum^\infty_{j_{p_e} = 0} c_{j_1, j_2, \ldots, j_m} \leqno(1.4)$$
is called a {\it subseries} of the $m$-multiple series (1.1). We refer to [8, Examples 1 and 2], where examples are given 
in the case of double series. 

Next, we recall the notion of regular convergence for multiple series. This notion was introduced in [5] and 
called there as ``convergence in the restricted sense" (see also in [6]). Given any pair 
$$(k_1, k_2, \ldots, k_m), \quad (l_1, l_2, \ldots, l_m) \in \N^m, \quad 0\le k_p\le l_p, \quad p=1,2,\ldots, m;$$
we set 
$$s(k_1, l_1; k_2, l_2; \ldots; k_m, l_m):= \sum^{l_1}_{j_1=k_1} \ \sum^{l_2}_{j_2=k_2} 
\ldots \sum^{l_m}_{j_m=k_m} c_{j_1, j_2, \ldots, j_m},\leqno(1.5)$$
which may be called {\it subrectangular sums} of (1.1). In particular, (1.5) reduces to the rectangular partial 
sum if $k_1 = k_2=\ldots = k_m=0$, that is, we have 
$$s(l_1, l_2, \ldots, l_m) = s (0, l_1; 0, l_2; \ldots; 0, l_m).$$

Now, the multiple series (1.1) is said to {\it converge regularly} if for every $\e>0$ there exists 
$\lambda_2 = \lambda_2 (\e) \in \N$ such that 
$$|s(k_1, l_1; k_2, l_2; \ldots; k_m, l_m)|<\e\leqno(1.6)$$
$${\rm if}\quad  \max\{k_1, k_2, \ldots, k_m\} > \lambda_2\quad {\rm and}\quad  0\le k_p \le l_p,   \ \ p=1,2,\ldots, m.$$
It is easy to check that if (1.6) is satisfied, then we have 
$$|s(k_1, k_2, \ldots, k_m) - s(l_1, l_2, \ldots, l_m)| < m\e.$$
Since $\e>0$ is arbitrary, by the Cauchy convergence criterion, it follows that the regular convergence of (1.1) 
implies its convergence in Pringsheim's sense. Consequently, the sum of a regularly convergent multiple series is well defined. On the other hand, a multiple series may converge in 
Pringsheim's sense without converging regularly. See, e.g., [8, Example 3] in the case of double series. 

It is obvious that if the multiple series (1.1) {\it converges absolutely}, that is, if 
$$\sum^\infty_{j_1=0} \ \sum^\infty_{j_2 = 0} \ldots \sum^\infty_{j_m=0} |c_{j_1, j_2, \ldots, j_m}| < \infty,$$
then it also converges regularly. The converse statement is not true in general. For example see 
[8, Example 5] in the case of double series. 

\bigskip 

\noindent {\bf Remark 1.2.} For double series, Hardy [1] introduced the notion of regular convergence as follows. 
The double series 
$$\sum^\infty_{j_1 = 0} \ \sum^\infty_{j_2 = 0} c_{j_1, j_2}$$
is said to converge regularly if it converges in Pringsheim's sense, and if each of its so-called row and column subseries defined by 
$$\sum^\infty_{j_2 = 0} c_{j_1, j_2}, \quad {\rm where} \quad j_1\in \N; \ \sum^\infty_{j_1=0} 
c_{j_1, j_2}, \quad {\rm where} \quad j_2\in \N;$$ 
converges as a single series. It is easy to check (see [6]) that Hardy's definition is equivalent to the definition given in 
(1.6) for $m=2$. 

For $m$-multiple series, where $m\ge 3$, and equivalent definition of regular convergence is formulated in the following 

\bigskip 

\noindent {\bf Theorem 1.1.} {\it The $m$-multiple series} (1.1) {\it converges regularly if and only if 
it converges in Pringsheim's sense and if each of its $(m-1)$-
multiple subseries 
$$\sum^\infty_{j_1=0} \ldots \sum^\infty_{j_{p-1} = 0} \sum^\infty_{j_{p+1} = 0} \ldots \sum^\infty_{j_m = 0} c_{j_1, j_2, \ldots, j_m}$$
converges regularly for all choices of $p\in \{1,2,\ldots, m\}$ and} $j_p\in \N$. 

We note that Theorem 1.1 was stated in [5, Theorem 1$'$] in a wrong form, while the correct statement appeared in [6, Theorem 1] without a detailed proof. Therefore, we sketch a 
proof below. 

\bigskip

\noindent {\bf Proof of Theorem 1.1.} {\it Necessity}. It is trivial. 

\noindent {\it Sufficiency.} It follows from (1.3) that given $\e>0$, we have 
$$|s(k_1, l_1; k_2, l_2; \ldots; k_m, l_m)|\leqno(1.7)$$
$$=\Big|\sum^1_{\delta_1=0} \ldots \sum^1_{\delta_m = 0} (-1)^{\delta_1+\ldots+\delta_m} s(\delta_1(k_1-1) + 1-\delta_1)l_1,$$
$$\ldots, \delta_m (k_m-1) + 1-\delta_m) \ell_m)| < 2^m\e$$
$${\rm if} \quad \min\{k_1, \ldots, k_m\} > \lambda_1 (e) \quad {\rm and} \quad 
0\le k_p \le l_p, \ p=1,2,\ldots, m;$$
with the agreement that 
$$s(k_1, k_2, \ldots, k_m):=0 \quad {\rm if} \quad \min\{k_1-1, \ldots, k_m-1\} = 0.\leqno(1.8)$$

For the sake of brevity in writing, we consider triple series and its double subseries 
$$\sum^\infty_{j_2=0} \ \sum^\infty_{j_3=0} c_{j_1, j_2, j_3}, \quad {\rm where} \quad j_1\in \{0,1,\ldots, \lambda_1\},$$
$$\sum^\infty_{j_1=0} \ \sum^\infty_{j_3=0} c_{j_1, j_2, j_3}, \quad {\rm where} \quad j_2\in \{0,1,\ldots, \lambda_1\},$$
$$\sum^\infty_{j_1=0} \ \sum^\infty_{j_2=0} c_{j_1, j_2, j_3}, \quad {\rm where} \quad j_3\in \{0,1,\ldots, \lambda_1\}.$$
By assumption, each of these double series converges regularly and the number of them is $(1+\lambda_1)^3$. Since 
a finite number of regularly convergent double series converge uniformly, there exists some $\lambda_3 = \lambda_3 (\e) \in \N$ 
such that 
$$|s(k_1, l_1; k_2, l_2; k_3, l_3)|<\e$$
$${\rm if} \quad \max\{k_2, k_3\}>\lambda_3 \quad {\rm and} \quad 0\le k_1\le l_1\le \lambda_1;$$
and two more analogous inequalities hold true. Combining these inequalities with 
(1.7) when $m=3$ yields 
$$|s(k_1, l_1; k_2, l_2; k_3, l_3)| < 10\e$$
$${\rm if} \quad \max\{k_1, k_2, k_3\} > \max\{\lambda_1, \lambda_3\} \quad {\rm and} \quad 
0\le k_p \le \ell_p,\quad  p=1,2,3.$$
Since $\e>0$ is arbitrary, this proves the regular convergence of the $m$-multiple series (1.1) in case $m=3$. $\kocka$ 

The next corollary immediately follows from the Sufficiency part of Theorem 1.1. 

\bigskip 

\noindent {\bf Corollary 1.1.} {\it The $m$-multiple series} (1.1) {\it converges regularly if and only if it 
converges in 
Pringsheim's sense and} 

(i) {\it each of its subseries defined in} (1.4) {\it converges regularly; or equivalently} 

(ii) {\it each of its subseries defined in} (1.4) {\it converges in Pringsheim's sense, 

\noindent
where by the regular convergence as well as by the convergence in Pringsheim's sense of a single series we mean 
its ordinary convergence. }

\bigskip

\noindent {\bf Remark 1.3.} It is of some interest to observe that both Theorem 1.1 and Corollary 1.1 remain valid if the convergence of the $m$-multiple series (1.1) 
in Pringsheim's sense is exchanged in them for the following weaker one: for every 
$\eta>0$ there exists some $\lambda_4 = \lambda_4 (\eta) \in \N$ such that condition (1.7) is satisfied with $\eta$ instead of $2^m \e$. 

\bigskip

\noindent {\bf Remark 1.4.} Similarly to the two convergence notions above for $m$-multiple series, we can define 
analogous convergence notions for multiple sequences of complex numbers. The $m$-{\it multiple sequence} 
$$(s_{l_1, l_2, \ldots, l_m} : (l_1, l_2, \ldots, l_m) \in \N^m) \subset \C$$
is said to {\it converge in Pringsheim's sense} to $s\in \C$ if for every 
$\e>0$ there exists $\kappa_1=\kappa_1 (\e) \in \N$ such that 
condition (1.3) is satisfied. Next, motivated by (1.5)-(1.8), the $m$-multiple sequence 
$(s_{l_1, l_2, \ldots, l_m})$ is said to {\it converge regularly} if for every $\e>0$ 
there exists $\kappa_2 = \kappa_2 (\e) \in \N$ such that 
$$\Big|\sum^1_{\delta_1=0} \ldots \sum^1_{\delta_m=0} (-1)^{\delta_1+\ldots+ \delta_m} s(\delta_1 (k_1-1) + 
(1-\delta_1) \ell_1, \ldots, \delta_m (k_m-1) + (1-\delta_m) l_m)\Big|<\e$$
$${\rm if} \quad \max\{k_1, \ldots, k_m\} > \kappa_2 \quad {\rm and} \quad 0\le k_p \le l_p, \quad 
p=1,2,\ldots, m.$$

It is routine to check that the regular convergence of an $m$-multiple sequence implies its convergence in 
Pringsheim's sense. Consequently, the finite limit of a regularly convergent multiple sequence is well defined. 
Furthermore, if the multiple sequence $(s(l_1, l_2, \ldots, l_m)$) converges regularly, then each of its subsequences 
$$\Big(s(l_1, l_2, \ldots, l_m): (l_{p_1}, l_{p_2}, \ldots, l_{p_e}) \in \N^e,$$
$${\rm where}\quad (j_{p_e+1}, \ldots, j_{p_m}) \in \N^{m-e}\quad {\rm  is \ fixed \ arbitrarily}\Big) $$
converges regularly (as to the notation, cf. (1.4)). 

\bigskip

\noindent {\bf Remark 1.5.} In Harmonic Analysis (for example, multiple Fourier series, see in [13, Ch. XVII]) we frequently meet $m$-multiple 
series of the form 
$$\sum^\infty_{j_1=-\infty} \ \sum^\infty_{j_2=-\infty} \ldots \sum^\infty_{j_m = -\infty} 
c_{j_1, j_2, \ldots, j_m}.\leqno(1.9)$$
Using {\it symmetric rectangular partial sums} defined by 
$$s(l_1, l_2, \ldots, l_m) : = \sum^{l_1}_{j_1 = -l_1} \ \sum^{l_2}_{j_2=-l_2} \ldots 
\sum^{l_m}_{j_m=-l_m} c_{j_1, j_2, \ldots, j_m},$$ 
where $(l_1, l_2, \ldots, l_m) \in \N^m$ convergence of the multiple series (1.9) in Pringsheim's sense is also defined by 
(1.3), and denoted by (1.2). 

In the definition of the regular convergence of (1.9), instead of (1.6), we require the fulfillment of the following condition 
(cf. notation (1.5)): 
$$\Big|\sum_{k_1\le |j_1|\le l_1} \ \sum_{k_2\le |j_2| \le l_2} \ldots 
\sum_{k_m \le |j_m|\le l_m} c_{j_1, j_2, \ldots, j_m}|<\e$$
$${\rm if} \quad \max\{k_1, k_2, \ldots, k_m\}>\lambda_2(\e) \quad {\rm and} \quad 
0\le k_p \le l_p, \quad p=1,2,\ldots, m.$$

\bigskip

\noindent {\bf Remark 1.6.} A third convergence notion of multiple series was introduced in [4, p. 34]. The $m$-multiple series (1.1) is said to {\it converge 
completely} if it converges in Pringsheim's sense and if ``every single series obtained from 
it by holding all the subscripts of the terms but one 
fixed, is convergent". 
Clearly, the notion of complete convergence and that of regular convergence coincide for double series. 
On the other hand, making use of Corollary 1.1, it is easy to construct a triple series 
that converges completely, but not regularly. 

\vglue1cm 

\noindent {\bf 2. New result: successive summation of regularly convergent multiple series} 

It is clear that if the multiple series (1.1) converges absolutely, then its sum can be also 
computed by successive summation; 
that is, for any permutation $\{\sigma(1), \sigma(2), \ldots, \sigma(m)\}$ of $\{1,2,\ldots, m\}$, we have 
$$\sum^\infty_{j_{\sigma(1)} = 0 } \ \Big(\sum^\infty_{j_{\sigma(2)} = 0} \Big(\ldots \Big(\sum^\infty_{j_{\sigma(m)}=0 }
c_{j_1, j_2, \ldots, k_m}\Big) \ldots \Big) \Big) = s,\leqno(2.1)$$
where $s$ is the sum of (1.1) in Pringsheim's sense. 

Now, our new result is the following 

\bigskip

\noindent {\bf Theorem 2.1.} {\it If the $m$-multiple series} (1.1) {\it converges regularly, then} 
(2.1) {\it holds true.}

{\it Proof.} We will prove (2.1) by induction on $m$. In the case of  double series, Theorem 2.1 
was proved in [8, 
Theorem 1]. As induction hypothesis, we assume that (2.1) holds true for some 
$m\ge 2$, and we will prove that then (2.1) also holds true for $m+1$. 
For the sake of brevity in writing, we present this induction step in the case where $m=2$ 
and $\sigma(p) = p$ for $p=1,2,3$. 

By Corollary 1.1, each subseries of (1.1) defined by (1.4) also converges regularly. 
In particular, for fixed $(j_1, j_2) \in \N^2$, denote by $d_{j_1, j_2}$ the sum of the single subseries 
$$\sum^\infty_{j_3=0} c_{j_1, j_2, j_3} = : d_{j_1, j_2}.\leqno(2.2)$$
We claim that the double series
$$\sum^\infty_{j_1=0} \ \sum^\infty_{j_2 = 0} d_{j_1, j_2}\leqno(2.3)$$
converges regularly. Indeed, it follows from (1.6) (cf. notation (1.5) that 
$$\Big|\sum^{l_1}_{j_1=k_1} \ \sum^{l_2}_{j_2=k_2} \ \sum^{l_3}_{j_3=0} c_{j_1, j_2, j_3}\Big| < \e 
\quad {\rm if} \quad \max\{k_1, k_2\} > \lambda_2 (\e),$$
$$0\le k_p \le l_p, \quad p=1,2,\quad {\rm and} \quad l_3\ge 0.$$
Letting $\ell_3 \to \infty$ and keeping (2.2) in mind gives 
$$\Big|\sum^{l_1}_{j_1=k_1} \sum^{l_2}_{j_2=k_2} d_{j_1, j_2}\Big|\le \e\quad {\rm if} \quad 
\max\{k_1, k_2\} > \lambda_2 (\e)\leqno(2.4)$$
$${\rm and} \quad 0\le k_p \le \ell_p, \quad p=1,2.$$
Since $\e>0$ is arbitrary in (2.4), we conclude that the double series (2.3) converges regularly. 

It follows also from (1.6) that 
$$\Big|\sum^{l_1}_{j_1=0} \ \sum^{l_2}_{j_2=0} \ \sum^{l_3}_{j_3=k_3+1} 
c_{j_1, j_2, j_3}\Big| < \e\leqno(2.5)$$
$${\rm if} \quad l_1, l_2\ge 0\quad {\rm and} \quad l_3>k_3 \ge \lambda_2 (\e).$$

Denote by $s$ the sum of the triple series 
$$\sum^\infty_{j_1=0} \ \sum^\infty_{j_2=0} \ \sum^\infty_{j_3=0} c_{j_1, j_2, j_3} = :s \leqno(2.6)$$
in Pringsheim's sense. By (1.3) and (2.4)-(2.5) we obtain that 
$$\Big|\sum^{l_1}_{j_1=0} \ \sum^{l_2}_{j_2=0} d_{j_1, j_2} -s\Big|$$
$$\le \Big|\sum^{l_1}_{j_1=0} \ \sum^{l_2}_{j_2=0} \ \sum^{k_3}_{j_3 =0} c_{j_1, j_2, j_3} - s\Big| 
+ \Big|\sum^{l_1}_{j_1=0} \ \sum^{l_2}_{j_2=0} \ \sum^\infty_{j_3=k_3+1} 
c_{j_1, j_2, j_3}\Big|$$
$$\le |s_{l_1, l_2, k_3} - s| + \e < 2\e\quad {\rm if} \quad \min\{l_1, l_2, k_3\} 
> \lambda_1 (\e) \quad {\rm and} \quad k_3\ge \lambda_2 (\e).$$
Since $\e>0$ is arbitrary, this proves that the double series (2.3) converges to $s$ 
 in Pringsheim's sense; or we may equivalently write that 
$$\sum^\infty_{j_1=0} \ \sum^\infty_{j_2=0} \Big(\sum^\infty_{j_3=0} c_{j_1, j_2, j_3}\Big) = 
\sum^\infty_{j_1=0} \ \sum^\infty_{j_2=0} \ \sum^\infty_{j_3=0} c_{j_1, j_2, j_3}.\leqno(2.7)$$ 

We have proved above that the double series (2.3) converges regularly (see (2.4)). Thus, we may apply the induction hypothesis to obtain that 
$$\sum^\infty_{j_1=0} \ \sum^\infty_{j_2=0} d_{j_2, j_2} = \sum^\infty_{j_1=0} 
\Big(\sum^\infty_{j_2=0} d_{j_1, j_2}\Big).\leqno(2.8)$$
Combining (2.7) and (2.8), while keeping (2.2) and (2.6) in mind, we conclude that 
$$\sum^\infty_{j_1=0} \ \sum^\infty_{j_2=0} \ \sum^\infty_{j_3=0} c_{j_1, j_2, j_3} = 
\sum^\infty_{j_1=0} \ \sum^\infty_{j_2=0} d_{j_1, j_2}$$
$$=\sum^\infty_{j_1=0} \Big(\sum^\infty_{j_2=0} d_{j_1, j_2}\Big) = 
\sum^\infty_{j_1=0} \Big(\sum^\infty_{j_2=0} \Big(\sum^\infty_{j_3=0} c_{j_1, j_2, j_3}\Big)\Big).$$
This proves (2.1) for regularly convergent triple series. 

The proof of (2.1) for regularly convergent $m$-multiple series with $m\ge 4$ goes along analogous lines by induction argument. 

The proof of Theorem 2.1 is complete. $\kocka$ 

\vglue1cm

\noindent {\bf 3. Convergence of multiple integrals of locally integrable functions} 

Let $f: \overline{\R}^m_+ \to \C$ be a locally integrable function 
in Lebesgue's sense over the closed positive octant $\overline{\R}^m_+ := [0, \infty)^m$, in symbols: 
$f\in L^1_{\loc} (\overline{\R}^m_+)$, where $m\ge 2$ is a fixed integer. We consider the 
$m$-{\it multiple integral} 
$$\int^\infty_0 \int^\infty_0 \ldots \int^\infty_0 f(t_1, t_2, \ldots, t_m) dt_1 dt_2 \ldots dt_m,\leqno(3.1)$$
whose $m$-fold {\it rectangular partial integrals} are defined by 
$$I(v_1, v_2, \ldots, v_m) := \int^{v_1}_0 \int^{v_2}_0 \ldots \int^{v_m}_0 
f(t_1, t_2, \ldots, t_m) dt_1 dt_2 \ldots dt_m,\leqno(3.2)$$
$${\rm where} \quad (v_1, v_2, \ldots, v_m) \in \overline{\R}^m_+$$

Analogously to the convergence of multiple series (cf. (1.2) and (1.3)), the multiple integral (3.1) is said to 
{\it converge in Pringsheim's sense} to the finite limit $I\in \C$, or equivalently, it is said that this $I$ is the {\it value} 
(or sum) of the multiple integral (3.1) in Pringsheim's sense, in symbols: 
$$\lim_{v_1, v_2, \ldots, v_m \to \infty} I(v_1, v_2, \ldots, v_m) = I,\leqno(3.3)$$
if for every $\e>0$ there exists $\rho_1 = \rho_1 (\e)\in \R_+$ such that 
$$|I(v_1, v_2, \ldots, v_m) - I| < \e \quad {\rm if} \quad 
\min\{v_1, v_2, \ldots, v_m\} > \rho_1.\leqno(3.4)$$

Next, we introduce the notion of regular convergence for $m$-multiple integrals. Given any pair 
$$(u_1, u_2, \ldots, u_m), \quad (v_1, v_2, \ldots, v_m) \in \overline{\R}^m_+, \quad 
0\le u_k \le v_k, \ \ k=1,2,\ldots, m;$$
we set
$$I(u_1, v_1; u_2, v_2; \ldots; u_m, v_m)\leqno(3.5)$$
$$:= \int^{v_1}_{u_1} \int^{v_2}_{u_2} \ldots \int^{v_m}_{u_m} f(t_1, t_2, \ldots, t_m) dt_1 dt_2 \ldots dt_m, $$
which may be called {\it subrectangular integrals} of (3.1). 

We note that the notion of regular convergence of double integrals was introduced in [7] by the present author. 

Now, the multiple integral (3.1) is said to {\it converge regularly} if for every $\e>0$ there 
exists $\rho_2 = \rho_2 (\e) \in \R_+$ such that 
$$|I(u_1, v_1; u_2 v_2; \ldots; u_m, v_m)|<\e\leqno(3.6)$$
$${\rm if} \quad \max\{u_1, u_2, \ldots, u_m\} > \rho_2 \quad {\rm and} \quad 0\le u_k \le v_k, \ 
k=1,2,\ldots, m.$$
In particular, (3.5) reduces to (3.2) if $u_k=0, k=1,2,\ldots, m;$ that is, we have 
$$I(v_1, v_2, \ldots, v_m) = I(0, v_1; 0, v_2; \ldots; 0, v_m).$$

It is easy to check that if (3.6) is satisfied, then 
$$|I(u_1, u_2, \ldots, u_m) - I(v_1, v_2, \ldots, v_m)| < m\e.$$
Since $\e>0$ is arbitrary, the application of  the Cauchy convergence criterion gives  that the regular convergence of (3.1) 
implies its convergence in 
Pringsheim's sense as indicated in symbols by (3.3). 
Consequently, the value of a regularly convergent multiple integral is well defined. On the other hand, a multiple integral 
may converge in Pringsheim's sense without converging regularly. 
See, e.g., [8, just before Example 6] in the case of a double integral. 

\bigskip 

\noindent {\bf Remark 3.1.} It is obvious that if $f\in L^1 (\overline{\R}^m_+)$ in (3.1), then the multiple integral 
(3.1) converges regularly and its value equals the Lebesgue integral of $f$ over the whole octant $\overline{\R}^m_+$. 
The converse statement is not true in general. See, e.g., [8, Example 6] in the case of a double integral. 

\bigskip

\noindent {\bf Remark 3.2.} Suppose that a function 
$$J:= J(u_1, v_1; u_2, v_2; \ldots; u_m, v_m) : \overline{\R}^{2m}_+\to \C, \quad 0\le u_k\le v_k, \ k=1,2,\ldots, m; \leqno(3.7)$$
is such that 
$$J=0 \quad {\rm if} \quad u_k = v_k \quad {\rm or} \quad v_k = 0\quad {\rm for\ some}\quad 1\le k\le m,\leqno(3.8)$$
and $J$ enjoys the property of {\it additivity} in each pair of its variables, by which we mean, e.g., in the case 
of $(u_1, v_1)$ the following: 
$$J(u_1, v_1; u_2, v_2; \ldots) = J(u_1, \widetilde u_1; u_2, v_2; \ldots) + J(\widetilde u_1, v_1; u_2, v_2; \ldots)\leqno(3.9)$$
$${\rm for\ all}\quad 0\le u_1 < \widetilde u_1 < v_1, \quad 0\le u_k\le v_k, \quad k=2,3,\ldots, m.$$
Under these conditions, by means of $J(u_1, v_1; u_2, v_2; \ldots; u_m, v_m)$, analogously to (3.5) we can define the regular convergence for the symbol 
$J(0, \infty; 0, \infty; \ldots; 0, \infty)$ (cf. (3.1); as well as by means of $J(0, v_1; 0, v_2; \ldots; 0, v_m)$, analogously to (3.4), we can also define 
its convergence in Pringsheim's sense to a finite limit. 

\bigskip 
\noindent {\bf Remark 3.3.} In Harmonic Analysis (for example, Fourier transform, see in [12, Ch.I]) we frequently meet multiple integrals 
of the form 
$$\int^\infty_{-\infty} \int^\infty_{-\infty} \ldots \int^\infty_{-\infty} f(t_1, t_2, \ldots, t_m) dt_1 dt_2\ldots dt_m, \leqno(3.10)$$
where $f\in L^1_{\loc} (\R^m)$. Using {\it symmetric rectangular partial integrals} (cf. (3.2)) defined by 
$$I(v_1, v_2, \ldots, v_m):= \int^{v_1}_{-v_1} \int^{v_2}_{-v_2} \ldots \int^{v_m}_{-v_m} f(t_1, t_2, \ldots, t_m) 
dt_1 dt_2\ldots dt_m,$$
$${\rm where}\quad (v_1, v_2, \ldots, v_m) \in \overline{\R}^m_+,$$
the convergence of (3.10) in Pringsheim's sense is defined also by (3.4). 

In the definition of regular convergence of the multiple integral (3.10), 
instead of (3.6) with the notation (3.5), we require the fulfillment of the following condition: for every 
$\e>0$ there exists $\rho_2 = \rho_2 (\e)\in \R_+$ such that 
$$\Big|\int_{u_1<|t_1|<v_1} \ldots \int_{u_m < |t_m| < v_m} f(t_1, \ldots, t_m) dt_1 \ldots dt_m\Big|<\e$$
$${\rm if} \quad \max\{u_1, \ldots, u_m\} > \rho_2 \quad {\rm and} \quad 0\le u_k\le v_k, \quad k=1,2,\ldots, m.$$

\vfill\eject

\noindent {\bf 4. New results: generalized versions of Fubini's theorem} 

We recall Fubini's classical theorem on the successive integration of an $m$-multiple integral (see [1,2] by 
Fubini, and see also, e.g., [11, p. 85] by F. Riesz and B. Sz.-Nagy), according to which if $f\in L^1 (\overline{\R}^m_+)$ and 
$m=p+q$, where $p,q\in \N_+$, then we have 
$$\int_{\overline{\R}^m_+} f(t_1, t_2, \ldots, t_m) dt_1 dt_2 \ldots dt_m\leqno(4.1)$$
$$=\int_{\overline{\R}^p_+} \Big(\int_{\overline{R}^q_+ } f(t_1, t_2, \ldots, t_m) dt_{p+1} dt_{p+2} \ldots dt_m\Big) dt_1 dt_2 \ldots dt_p,$$
where the inner integral exists in Lebesgue's sense for almost every $(t_1, t_2, \ldots, t_p) \in \overline{\R}^p_+$, 
and the outer integral exists in Lebesgue's sense. 

\bigskip 

\noindent {\bf Remark 4.1.} Fubini's theorem holds even in the following more general form:If $f\in L^1(\overline{\R}^m_+)$ and $\{\sigma(1), \sigma(2), \ldots, 
\sigma(m)\}$ is a permutation of $\{1,2,\ldots, m\}$, then the left-hand side in (4.1) can be computed by 
the following successive integration: 
$$\int^\infty_0 \Big(\int^\infty_0 \Big(\ldots \Big(\int^\infty_0 f(t_1, t_2, \ldots, t_m) dt_{\sigma(m)}\Big) 
\ldots\Big) dt_{\sigma(2)}\Big) dt_{\sigma(1)}.\leqno(4.2)$$

Our goal is to prove a generalized version of (4.1) under the weaker assumptions that $f\in L^1_{\loc} (\overline{\R}^m_+)$ 
and that the $m$-multiple integral on the left-hand side of (4.1) converges regularly. 

\bigskip

\noindent {\bf Remark 4.2.} If $f\in L^1_{\loc} (\overline{\R}^m_+)$, then by definition $f\in L^1{\cal R}_m)$ for every $m$-fold bounded rectangle 
$${\cal R}_m:= [a_1, b_1] \times [a_2, b_2]\times \ldots \times [a_m, b_m], \quad 
0\le a_k<b_k, \ \ k=1,2,\ldots, m.$$
Since any countable union of sets of Lebesgue measure zero in a Euclidean space is also of measure zero, by virtue of 
Fubini's theorem, for the function $f(t_1, \ldots, t_p, t_{p+1}, \ldots, t_m) \in L^1_{\loc} (\overline{\R}^q_+)$ for almost every 
$(t_1, \ldots, t_p)\in \overline{\R}^p_+$. 

Our main new result is the following

\bigskip

\noindent {\bf Theorem 4.1.} {\it If $f\in L^1_{\loc} (\overline{\R}^m_+)$, the $m$-multiple integral} (4.1) {\it converges 
regularly, and $m=p+q$, where $p, q\in \N_+$, then the finite limit 
$$\lim_{v_{p+1}, v_{p+2}, \ldots, v_m \to \infty} I(u_1, v_1; \ldots; u_p, v_p; 0, v_{p+1}; 0, v_{p+2}; \ldots; 0, v_m)\leqno(4.3)$$
$$=: J(u_1, v_1; u_2, v_2; \ldots; u_p, v_p), \quad 0\le u_k\le v_k, k=1,2,\ldots, p;$$
exists uniformly in each of its variables. Furthermore, the function $J: \overline{R}^{2p}_+ \to \C$ enjoys the property of additivity 
as indicated by} (3.7) - (3.9) {\it in Remark} 3.2, {\it converges regularly in the sense that for every 
$\e>0$ there exists some $\rho_3 = \rho_3(\e) \in \R_+$ such that 
$$|J(u_1, v_1; u_2, v_2; \ldots; u_p, v_p)|\le \e\leqno(4.4)$$
$$if\quad \max\{u_1, u_2, \ldots, u_p\}>\rho_3 \quad and\quad 0\le u_k \le v_k, \quad k=1,2,\ldots, p;$$
and the finite limit 
$$\lim_{v_1, v_2, \ldots, v_p\to \infty} J(0, v_1; 0, v_2; \ldots; 0, v_p)=I\leqno(4.5)$$
also exists, where $I$ is the finite limit of the multiple integral} (3.1) {\it in Pringsheim's sense. } 
{\bf Proof.} For the sake of brevity in writing, we present the proof in the case when $m=4$ and $p=q=2$. That is, from now on, we consider the 
{\it quadruple integral} 
$$\int^\infty_0 \int^\infty_0\int^\infty_0\int^\infty_0 f(t_1, t_2, t_3, t_4) dt_1 dt_2 dt_2 dt_4,\leqno(4.6)$$
which by assumption converges regularly. By definition, for every $\e>0$ there exists $\rho_2=\rho_2(\e)\in \R_+$ such that 
$$\Big|\int^{v_1}_{u_1} \int^{v_2}_{u_2} \int^{v_3}_{u_3} \int^{v_4}_{u_4} f(t_1, t_2, t_3, t_4) dt_1 dt_2 dt_3 dt_4\Big|<\e\leqno(4.7)$$
$${\rm if}\quad \max\{u_1, u_2, u_3, u_4\}> \rho_2 \quad {\rm and} \quad 0\le u_k \le v_k, \quad k=1,2,3,4.$$
For fixed $0\le u_k \le v_k$, $k=1,2$; this means that 
$$\Big|\int^{v_1}_{u_1} \int^{v_2}_{u_2} \Big(\int^{v_3}_{u_3} \int^{v_4}_{u_4} \ f(t_1, t_2, t_3, t_4) 
dt_3 dt_4) dt_1 dt_2\Big|< \e$$
$${\rm if} \quad \max\{u_3, u_4\}>\rho_2 \quad {\rm and} \quad 0\le u_k\le v_k, \quad k=3,4.$$
Hence it follows that 
$$\Big|\int^{v_1}_{u_1} \int^{v_2}_{u_2}\Big(\Big\{\int^{v_3}_0 \int^{v_4}_0 - \int^{\widetilde v_3}_0 \int^{\widetilde v_4}_0\Big\} f(t_1, 
t_2, t_3, t_4) dt_3 dt_4\Big)dt_1 dt_2\Big|<2\e$$
$${\rm if} \quad \min\{v_3, v_4, \widetilde v_3, \widetilde v_4\} > \rho_2,$$
independently of the values of $u_1, v_1; u_2, v_2$. Since $\e>0$ is arbitrary, by the Cauchy convergence criterion, we conclude that the (inner) 
double integral 
$$\int^{v_1}_{u_1} \int^{v_2}_{u_2}\Big(\int^\infty_0 \int^\infty_0 f(t_1, t_2, t_3, t_4) dt_3 dt_4\Big) dt_1 dt_2$$ 
converges in Pringsheim's sense, and even uniformly in 
$u_1, v_1; u_2, v_2$. That is, the finite limit 
$$\lim_{v_3, v_4\to \infty} \int^{v_1}_{u_1} \int^{v_2}_{u_2}\Big(\int^{v_3}_0 \int^{v_4}_0 
f(t_1, t_2, t_3, t_4) dt_3 dt_4=: J(u_1, v_1; u_2, v_2)$$
exists uniformly for all $0\le u_k \le v_k$, $k=1,2$. This proves (4.3) in the case, where $m=4$ and $p=q=2$. 

It follows from the additivity property of an integral that the limit function
$J(u_1, v_1; u_2, v_2)$ also enjoys the additivity property indicated in Remark 3.2.

Next, we claim that the limit function $J(u_1, v_1; u_2, v_2)$ converges regularly 
(in the sense as  indicated in Remark 3.2); or in an equivalent formulation, we may say that the 
outer double integral 
$$\int^\infty_0 \int^\infty_0 \Big(\int^\infty_0 \int^\infty_0 f(t_1, t_2, t_3, t_4) dt_3 dt_4\Big) dt_1 dt_2$$
converges regularly. To justify this claim, we start with (4.8), according to which for every 
$\e>0$ there exists some $\rho_4=\rho_4(\e)\in \R_+$ such that 
$$\Big|\int^{v_1}_{u_1} \int^{v_2}_{u_2} \Big(\int^{v_3}_0 \int^{v_4}_0 f(t_1, t_2, t_3, t_4) dt_3 dt_4\Big) dt_1 dt_2\leqno(4.9)$$
$$-J(u_1, v_1; u_2, v_2)|<\e\quad {\rm if} \quad \min\{v_3, v_4\}>\rho_4 \quad {\rm and} \quad 0\le u_k\le v_k, \ k=1,2.$$
Combining (4.7) and (4.9) gives
$$|J(u_1, v_1; u_2, v_2)|\leqno(4.10)$$
$$\le \Big|J(u_1, v_1; u_2, v_2) - \int^{v_1}_{u_1} \int^{v_2}_{u_2} 
\Big(\int^{v_3}_0 \int^{v_4}_0 f(t_1, t_2, t_3, t_4) dt_3 dt_4\Big) dt_1 dt_2\Big|$$
$$+\Big|\int^{v_1}_{u_1} \int^{v_2}_{u_2} \Big(\int^{v_3}_0 \int^{v_4}_0 f(t_1, t_2, t_3, t_4) dt_3 dt_4\Big) dt_1 dt_2\Big|<2\e$$
$${\rm if} \quad \max\{u_1, u_2\}>\rho_2\quad {\rm and} \quad \min\{v_3, v_4\}>\rho_4.$$
Since $\e>0$ is arbitrary, the fulfillment of (4.4) with $2\e$ in place of $\e$ follows from (4.10) in the case, where $m=4$ and 
$p=q=2$. 

We will denote by $I$ the finite limit of the quadruple integral (4.6) in Pringsheim's sense, which certainly exists, 
due to the assumption that (4.6) converges regularly. Thus, for every $\e>0$ there exists $\rho_1=\rho_1(\e)\in\R_+$ such that 
$$\Big|\int^{v_1}_0 \int^{v_2}_0 \int^{v_3}_0 \int^{v_4}_0 f(t_1, t_2, t_3, t_4) dt_1 dt_2 dt_3 dt_4-I\Big|<\e\leqno(4.11)$$
$${\rm if} \quad \min\{v_1, v_2, v_3, v_4\} > \rho_1.$$

Next, we fix $0\le u_k\le v_k$ for $k=1,2$; and take auxiliary variables $\widetilde v_k (>v_k)$ for $k=3,4$. Since 
$$\int^{v_1}_{u_1} \int^{v_2}_{u_2} \Big(\int^{v_3}_0 \int^{v_4}_0 f(t_1, t_2, t_3, t_4) dt_3 dt_4\Big) dt_1 dt_2$$
$$=\int^{v_1}_{u_1} \int^{v_2}_{u_2} \Big(\Big\{\int^{\widetilde v_3}_0 \int^{\widetilde v_4}_0 - 
\int^{\widetilde v_3}_{v_3} \int^{\widetilde v_4}_0 - \int^{v_3}_0 \int^{\widetilde v_4}_{v_4} 
\Big\} f(t_1, t_2, t_3, t_4)dt_3 dt_4\Big) dt_1 dt_2,$$
we conclude from (4.7) and (4.10) that 
$$\Big|\int^{v_1}_{u_1} \int^{v_2}_{u_2} \Big(\int^{v_3}_0 \int^{v_4}_0 f(t_1, t_2, t_3, t_4) dt_3 dt_4\Big) 
dt_1 dt_2 - J(u_1, v_1; u_2, v_2)\Big|\leqno(4.12)$$
$$\le \Big|\int^{v_1}_{u_1} \int^{v_2}_{u_2} \Big(\int^{\widetilde v_3}_0 \int^{\widetilde v_4}_0 
f(t_1, t_2, t_3, t_4) dt_3 dt_4\Big) dt_1 dt_2 - J (u_1, v_1; u_2, v_2)\Big|$$
$$+\Big|\int^{v_1}_{u_1} \int^{v_2}_{u_2} \Big(\int^{\widetilde v_3}_{v_3} \int^{\widetilde v_4}_0 f(t_1, t_2, t_3, t_4) dt_3 dt_4\Big) 
dt_1 dt_2\Big|$$
$$+\Big|\int^{v_1}_{u_1} \int^{v_2}_{u_2} \Big(\int^{v_3}_0 \int^{\widetilde v_4}_{v_4} f(t_1, t_2, t_3, t_4) dt_3 dt_4\Big) 
dt_1 dt_2\Big|< 4\e$$
$${\rm if} \quad \min\{v_3, v_4\} \ge \rho_2 \quad {\rm and} \quad \min\{\widetilde v_3, \widetilde v_4\}>\rho_4.$$
Combining (4.11) with (4.12) (we observe that $\widetilde v_3, \widetilde v_4$ are dummy variables in it), the latter one in the special case 
where $u_1=u_2 =0$ yields 
$$|J(0, v_1; 0, v_2) - I|\leqno(4.13)$$
$$\le \Big|J(0, v_1; 0, v_2) - \int^{v_1}_0 \int^{v_2}_0 \Big(\int^{v_3}_0 \int^{v_4}_0 f(t_1, t_2, t_3, t_4) dt_3 dt_4\Big) dt_1 dt_2\Big|$$
$$+\Big|\int^{v_1}_0 \int^{v_2}_0 \int^{v_3}_0 \int^{v_4}_0 f(t_1, t_2, t_3, t_4) dt_1 dt_2 dt_3 dt_4-I\Big| < 5\e$$
$${\rm if} \quad \min\{v_1, v_2\} > \rho_1 \quad {\rm and} \quad \min\{v_3, v_4\} > \max\{\rho_1, \rho_2\}.$$
Since $\e>0$ is arbitrary in (4.13), this proves (4.5) with $5\e$ in place of $\e$ in the case, where $m=4$ and $p=q=2$. 

The proof of Theorem 4.1 is complete. $\kocka$ 

A trivial corollary of Theorem 4.1 is the following 

\bigskip

\noindent {\bf Corollary 4.2.} {\it If $f\in L^1_{\loc} (\overline{\R}^m_+)$, the $m$-multiple integral} (3.1) 
{\it converges regularly, and 
$f(t_1, t_2, \ldots, t_p$, $t_{p+1}, \ldots, t_m)\in L^1 (\overline{\R}^q_+)$ for almost every 
$(t_1, t_2, \ldots, t_p) \in \overline{\R}^p_+$, where 
$m=p+q$ and $p,q\in \N_+$, then the $p$-multiple integral 
$$\lim_{v_{p+1}, v_{p+2}, \ldots, v_m} I(u_1, v_1; u_2, v_2; \ldots; u_p, v_p; 0, v_{p+1}; 0, v_{p+2}; \ldots; 0, v_m)\leqno(4.3')$$
$$=\int^{v_1}_{u_1} \int^{v_2}_{u_2} \ldots \int^{v_p}_{u_p} \Big(\int_{\overline{\R}^q_+} f(t_1, t_2, \ldots, t_p, 
t_{p+1}, \ldots, t_m) \times$$
$$\times dt_{p+1} dt_{p+2} \ldots dt_m\Big) dt_1 dt_2\ldots dt_p =: J(u_1, v_1; u_2, v_2; \ldots; u_p, v_p),$$
the function $J: \overline{\R}^{2p}_+ \to \C$ converges regularly, and the finite limit 
$$\lim_{v_1, v_2, \ldots, v_p} J(0, v_1; 0, v_2; \ldots; 0, v_p)=I\leqno(4.5')$$
also exists, where $I$ is the limit of the multiple integral} (3.1) {\it in Pringsheim's sense}. 

\bigskip 

\noindent {\bf Remark 4.3.} Under the conditions of Theorem 4.1, the function 
$J=J(u_1, v_1; u_2, v_2;$ $\ldots$; $u_p, v_p) : \overline{\R}^{2p}_+ \to \C$ 
enjoys the property of additivity, converges to $I$ in Pringsheim's sense (see (4.5')). Thus, a procedure analogous 
to the one described in Theorem 4.1 can be repeated for the function 
$J: \overline{\R}^{2p}_+ \to \C$ with $p=p_1+p_2$ and $p_1, p_2\in \N_+$. This procedure can be even repeated until we get to a function 
$J: \overline{\R}^2_+ \to \C$. But it is not clear to us how to interpret the whole process in terms of the 
traditional framework of integration. 

\bigskip

\noindent {\bf Remark 4.4.} If $p=m-1$ and $q=1$ in Theorem 4.1, we can repeat 
the process of successive integration for $J=J(u_1, v_1; 
u_2, v_2; \ldots; u_{m-1}, v_{m-1}): \overline{\R}^{2(m-1)}_+$ with choosing $p_1=m-2$ and $p_2=1$; and we repeat 
the process until we get to a counterpart of (4.2) with $\sigma(k) = k$, 
$k=1,2,\ldots, m$. But it is not clear to us again how to find 
an appropriate interpretation and notation for the whole process. 

\vglue1cm

\noindent {\bf References} 

\item{[1]} G. Fubini, Sugli integrali multipli, Atti Real Accad. Lincei Rend. Cl. Sci. Fis. 
Mat. Natur. (Roma) (5) 16(1907), no.1, 608-614. 
\item{[2]} G. Fubini, Sugli integrali doppii, ibid. (5) 22(1913), no. 1, 584-589. 
\item{[3]} G.H. Hardy, On the convergence of certain multiple series, Proc. Cambridge Philosoph. Soc., 
19(1916-1919), 86-95. 
\item{[4]} Ch. N. Moore, Summable series and convergence factors, Amer. Math. Soc., New York, 1938. 
\item{[5]} F. M\'oricz, On the convergence in a restricted sense of multiple series, 
Analysis Math., 5(1979), 135-147. 
\item{[6]} F. M\'oricz, Some remarks on the notion of regular convergence of 
multiple series, Acta Math. Hungar., 41 (1983), 161-168. 
\item{[7]} F. M\'oricz, On the uniform convergence of double sine integrals over $\overline{\R}^2_+$, 
Analysis, 31 (2011), 191-204. 
\item{[8]} F. M\'oricz, On the convergence of double integrals and a generalized version of Fubini's 
theorem on successive integration, Studia Math. (submitted for publication). 
\item{[9]} A. Pringsheim, Elementare Theorie der unendlichen Doppelreihen, M\"unch. Ber., 27 (1897), 101-132.
\item{[10]}  A. Pringsheim, Zur Theorie der zweifach unendlichen Zahlenfolgen, Math. Ann., 53(1900), 289-321.
\item{[11]} F. Riesz et B.Sz.-Nagy, Lecons d'analyse fonctionelle, Gauthier-Villars, Paris, 1955. 
\item{[12]} E.M. Stein and G. Weiss, Introduction to Fourier Analysis on Euclidean Spaces, Princeton Univ. Press, 1971. 
\item{[13]} A. Zygmund, Trigonometric Series, Vol. II, Cambridge Univ. Press, 1959. 

\bye